\documentclass{article}
\usepackage[T2A]{fontenc}
\usepackage[cp1251]{inputenc} 
\usepackage[english]{babel}
\usepackage{amsfonts,amssymb,mathrsfs,amscd}
\usepackage{amsthm,graphicx}

\newtheorem*{theorem}{Theorem}

\theoremstyle{definition}

\newtheorem*{st}{Statement}

\def\R{{\mathbb R}}

\title{Shear coordinates and braid invariants}

\author{Vassily Olegovich Manturov}

\begin{document}

\maketitle

\begin{abstract}
In the present paper we relate {\em shear coordinates} on
hyperbolic plane to braid invariants. 
\end{abstract}



\section{Introduction}

In \cite{InvariantsAndPictures}, when considering Voronoi  
diagrams dual to Delaunay triangulations, 
the author formalised the following 
\begin{st}
Octagon and far commutativity and yield braid invariants.
\end{st}

In \cite{MW,Roh}, a way of constructing braid invariants from solutions
to the octagon equation was established.

In that work, a braid represented by a dynamical system (a motion of $n$ points
in $\R^{2}$ or in $\R{}P^{2}$)
was split into ``generic'' parts separated by {\em flips}.
With each generic part we associate its Delaunay triangulation which does
not change as the picture changes generically and changes by a simple flip.

The main idea of \cite{InvariantsAndPictures}, see page 310, is:
to associate labels to edges of Delaunay triangulations
and whenever a triangulation undergoes a flip, 
apply the Ptolemy transformation to the corresponding diagonal
leaving other labels unchanged. It is well known that ``Ptolemy
transformation satisfies pentagon relation'',
as well as ``shear choordinate transformation satisfies pentagon
relation.''

\begin{figure}
\centering\includegraphics[width=200pt]{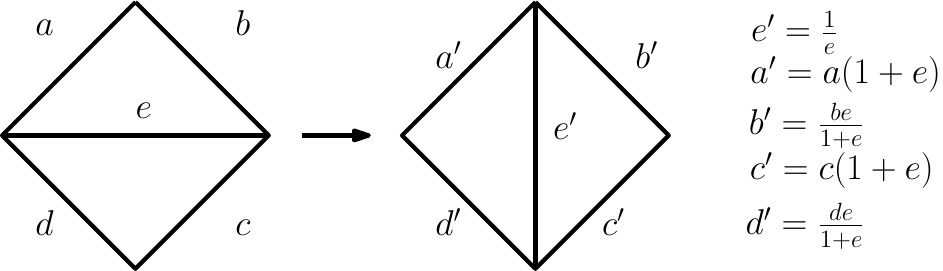}
\caption{The shear transformation of labels}
\label{shear}
\end{figure}

Unlike \cite{MW,Roh}, in the present paper the flips change coordinates
not only on the flipped edges but also on
the edges of the two triangles around it, see Fig. \ref{shear}.

Luckily, the operations of label changes commute if we perform flips inside
two neighbouring quadrilaterals, see Fig. \ref{Flipscommute} 
(if we take quadrilaterals which do not share any edge, the
flips commute by definition). 

\begin{frame}

\label{Flipscommute}
\end{frame}

\section{From braids to Voronoi 
 diagrams}

In the present section we represent an $n$-strand braid in $\R^{2}$ as
a dynamical system representing a motion of $n$ points.

We closely follow \cite{InvariantsAndPictures}.

Let $z_{i}(t),i=1,\cdots, n,t\in [0,1]$ be moving points on the plane $\R^{2}$.
For each $t$, we define the region
$U_{i}(t)$ to be $U_{i}(t)=\{z\in \R^{2}|\forall j:|z-z_{i}(t)|\le |z-z_{j}(t)|\}$.

Generically, these regions are separated by a trivalent graph
$\Gamma_{i}(t)$ with some infinite edges.
This graph is called the {\em Voronoi diagram} for $z(t)=\{z_{1}(t),
\dots, z_{n}(t)\}$.

The dual graph $D_{t}$ (which generically consists of
triangles) is called the {\em Delaunay triangulation}
(we do not pay attention to the ``infinite vertex'').

As points move, a typical {\em codimension 1} singularity
corresponds to those moments $t'$ for which some
four different $U_{j}(t')$ share a point.
This happens when some four {\em neighbouring} points from $z(t')$ belong
to the same circle (these points are neighbouring if
no other $z_{k}$ lies inside the circle).

Hence, a typical braid $\beta(t) = z(t)=\{z_{1}(t),\dots, z_{n}(t)\},
t\in [0,1]$ contains only finitely many moments $t'_{k}$ when
the Voronoi diagram is not generic.

We call such a braid {\em generic} if for each such $t'_{k}$
the diagram undergoes a flip, where
one one diagonal of a quadrilateral is
replaced with the other diagonal.

\subsection{Acknowledgements}

I am extremely grateful to Michael Shapiro for many fruitful
discussions concerning shear coordinates.
My special thanks to Platon Marulev and Daniil Tereshkin who helped me drawing pictures.

\section{The construction of the main invariant}

We fix some generic set of unordered $n$ points
to be $z'(0)\sim z'(1)$ for each braid to consider.

For this set  of points we have the Delaunay triangulation
$T=T(0)=T_{i}$.

We mark all edges of the triangulation by independent
variables $a_{1}, a_{2} \cdots $ (the number of variables
can be easily calculated from the Euler characteristic).

Whenever we undergo a flip, one edge is replaced by another
and in \cite{InvariantsAndPictures} the label of
the new edge is expressed in terms of 
labels of the former edges as follows

$$y=\frac{ac+bd}{x},$$

see Fig.\ref{Ptolemy}, and the other labels remain unchanged:

\begin{figure}
\centering\includegraphics[width=250pt]{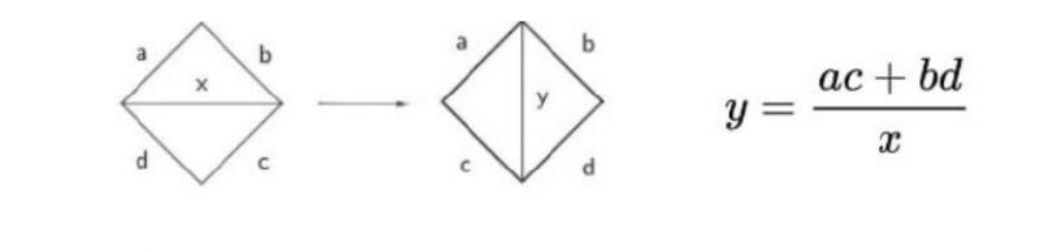}
\caption{The Ptolemy transformation of labels}
\label{Ptolemy}
\end{figure} 

We recall that the {\em shear transformation} of coordinates is shown in Fig. \ref{shear}.

Note that here not only the new edge gets a new label

$e'=\frac{1}{e}$ but also all four edges adjacent to it
change as follows:

$a'=a(1+e),b'=\frac{be}{1+e},c'=c(1+e), d'=\frac{de}{1+e}$.

Now, for a generic braid $\beta(t)$ we do the following:

For each flip we replace all labels as
above.

As the triangulations for $t=0$ and $t=1$ coincide,
we get a collection of new labels
$(a'_{1},a'_{2},\cdots )$.

Thus, we define

$T(\beta): (a_{1},a_{2}, \cdots, )\to
(a'_{1},a'_{2},\cdots)$.

Here $(a'_{1},a'_{2},\cdots)$ are some rational functions
in $(a_{1},a_{2} \cdots )$. In fact, from cluster algebra theory it is
known that they are Laurent polynomial in initial coordinates.

\section{The main theorem}

Collecting the above statements we get the following
\begin{theorem}
For two isotopic generic braids $\beta$ and $\beta'$
the transformations 
$T(\beta)$ and $T(\beta')$ are identical.
\end{theorem}

\begin{proof}

We consider an {\em  isotopy} $\beta$ and $\beta'$.

Denote this isotopy by 
$\beta_{s}, \beta_{0}=\beta,\beta_{1}=\beta'$.

An isotopy $\beta_{s}$ is {\em generic} if braids
$\beta_{s}$ are generic for all values of
$s$ except some $s_{1},\cdots, s_{p}$, and for
those values one of the codimension two events
happens:

\begin{enumerate}

\item the braid $\beta_{s_j}$ is typical but not
generic for some $t$ so that $\beta_{s_j- \varepsilon}$
and $\beta_{s_j+\varepsilon}$ change by a 
{\em back and forth transformation}:
we start with a triangulation
$(A,B,C),(C,D,A)$, change it to
$(B,C,D),(D,A,B)$ and return back.

\item for some value $t$, the set $\beta_{s_j}(t)$
contains some five neighbouring points on the same circle,
so that $\beta_{s_j- \varepsilon)}$ and $\beta_{s_j+ \varepsilon})$
differ by the following {\em Pentagon transformation},
see Fig. \ref{Pentagon}.

\item for some value $t$ the set $\beta_{s_j}(t)$ 
the flip happens in two places, so that 
$\beta_{s_j- \varepsilon)}$ and $\beta_{s_j+ \varepsilon)}$
differ by the following {\em commutativity},
see Fig. \ref{Commutativity}.

 \begin{figure}
 \centering\includegraphics[width=200pt]{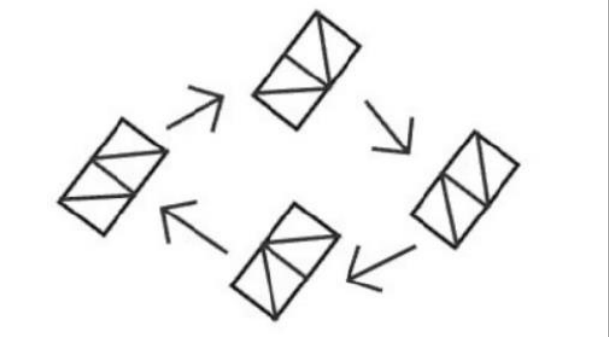}
   \caption{The far commutativity transformation} 
 \label{Commutativity}
 \end{figure}

\end{enumerate}

 \begin{figure}
 \centering\includegraphics[width=200pt]
 {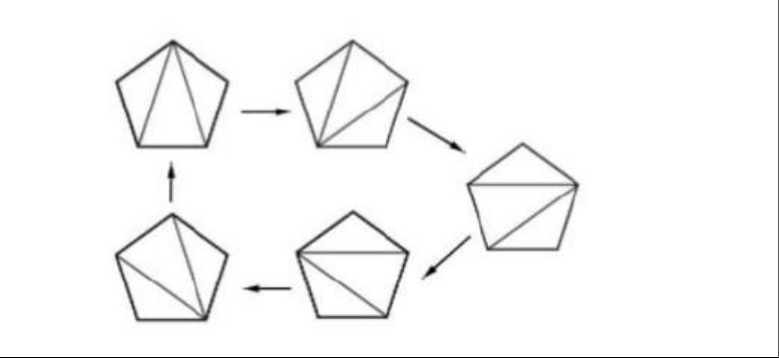}
   \caption{The pentagon transformation} 
 \label{Pentagon}
 \end{figure}

Now we claim that for each non-generic value 
$\beta_{s_j}$ the sequence of label transformation  for $\beta_{s_j}-\varepsilon$
and $\beta_{s_j}-\varepsilon$ in the neighbourhood of
value $t$ leads to the same result.

The case ``back and forth'' is obvious: on the left
hand side (say, $t-\varepsilon'$)
we the labels do not change at all
whence for $(t+\varepsilon')$ the labels change
twice and undergo two inverse transformations.

For the Ptolemy case, the {\em far commutativity} is obvious:
we have changes of two diagonals for two quadrilaterals
so that the edges of the quadrilaterals themselves do
not change, so the order of these two transformations
does not matter.

For the ``shear'' case the labels of edges of the quadrilateral
do change, however, it is known that
{\em shear coordinate transformations enjoy far commutativity}.

Finally, both Ptolemy transformation and
shear coordinate transformation satisfy the pentagon relation.

\end{proof}

\end{document}